\documentclass[12pt]{amsart}

\usepackage{vmargin}
\setmargrb{1in}{1in}{1in}{1in}
\usepackage{amstext}
\usepackage{amssymb,mathrsfs}
\usepackage{latexsym}
\newcommand\RR{{\mathbb R}}

\newcommand\NN{{\mathbb N}}
\newcommand\ZZ{{\mathbb Z}}

\newcommand{\cH}{\mathcal{H}}
\def\sp{\supseteq}

\def\inv{^{-1}}

\def\d={\,:=\,}

\newcommand{\semdir}
{\rtimes}

\font\frakten=eufm10
\newfam\frakfam
\textfont\frakfam=\frakten

\newcommand{\bbin}{\binom{2n}{n}^{1/2}}
\newcommand{\cU}{\mathcal{U}}
\newtheorem{thm}{Theorem}[section]
\newtheorem{lemma}[thm]{Lemma}
\newtheorem{cor}[thm]{Corollary}
\newtheorem{prop}[thm]{Proposition}
\newtheorem{Defn}[thm]{Definition}
\newtheorem{Ex}[thm]{Example}
\newtheorem{Rem}[thm]{Remark}
\newtheorem{Exs}[thm]{Examples}
\newtheorem{Rems}[thm]{Remarks}
\newtheorem{Defrem}[thm]{Definition and Remark}
\newtheorem{Remnt}[thm]{}
\newenvironment{defn}
 {\begin{Defn} \begin{rm}} {\end{rm} \hfill $\Box$ \end{Defn}}

\newenvironment{rem}
 {\begin{Rem} \begin{rm}} {\end{rm} \hfill $\Box$ \end{Rem}}

\newenvironment{prf} {{\bf Proof.}}{\hfill $\Box$}

\begin{document}

\author{H.~F{\"u}hr
\and K.~Gr\"ochenig}

\address{Institute of Biomathematics and Biometry \\
GSF National Research Center for Environment and Health \\
Ingolst{\"a}dter Stra{\ss}e~1 \\
D-85764 Neuherberg \\
Germany.}
\email{fuehr@gsf.de}
\email{karlheinz.groechenig@gsf.de}

\title[Sampling theorems from oscillation estimates]{Sampling theorems on locally compact groups from oscillation
estimates}

\keywords{Stratified Lie groups; Paley-Wiener space; sampling
theorems; frames; discrete series representations}
\subjclass[2000]{Primary 43A80, 42B35; Secondary 42C40, 26D10}

\date{\today}

\begin{abstract}
We present a general approach to derive sampling theorems on
locally compact groups from oscillation estimates.
We focus on the ${\rm
L}^2$-stability of the sampling operator by  using notions from frame
theory. This approach yields particularly simple and transparent
reconstruction procedures. We then apply these methods to the
discretization of discrete series representations  and to
Paley-Wiener spaces on stratified Lie groups.
\end{abstract}

\maketitle

\section{Introduction}  \label{sect:intro}

The sampling theorem of Shannon-Whittaker-Kotel'nikov is the prototype
of any sampling theorem. It states that  a function $f \in {\rm
L}^2(\RR)$ with $\mathrm{supp}\,\,\widehat{f}\, \subseteq [-1/2,1/2]$
 can be reconstructed from its sampled values $f|_\ZZ$ by
the cardinal series
\begin{equation} \label{eqn:shannon}
 f(x) = \sum_{k \in \ZZ} f(k) \frac{\sin \pi (x-k) }{ \pi (x-k)}~~,
\end{equation}
with convergence both in the ${\rm L}^2$-norm and uniformly. A related
 property is the norm equality
 \begin{equation} \label{eqn:shannon_norm} \| f
\|^2_2 = \sum_{k \in \ZZ} |f(k)|^2 ~~. \end{equation}

In communication theory and signal processing \eqref{eqn:shannon} is
considered the paradigm of a digital-analog conversion, in applied
mathematics the study of sampling and reconstruction theorems has
become a independent and very active field (as is documented by the
collections~\cite{ben01,ben04}). But the Shannon
sampling theorem is also important in investigations that are not
primarily
motivated by applications. Sampling and interpolation problems form a
whole branch of complex analysis~\cite{seip04}. Since by a theorem of
Paley and Wiener a function is bandlimited if and only if it possesses
an extension to an entire function of exponential type, spaces of
bandlimited functions are often called \emph{Paley-Wiener spaces}.

 In this paper, we adapt the point of view
that the sampling theorem is a phenomenon in  harmonic analysis; the
search for generalizations of the cardinal series in (nonabelian)
locally compact groups is an interesting  problem in itself,   and
abstract sampling theorems reveal new aspects of  the  analysis on
locally compact groups.

The Shannon-Wittaker-Kotel'nikov sampling theorem~\eqref{eqn:shannon}
and~\eqref{eqn:shannon_norm} follows directly from the Poisson summation
formula and is thus in the realm of locally compact abelian
groups. Indeed, a version of the sampling theorem for LCA groups was
formulated early on by Kluvanek~\cite{kluvanek}.

If $\RR $ is replaced by a nonabelian locally compact group $G$, one is
faced with  several fundamental questions:

(a) What is the appropriate concept of a bandlimited function on
$G$?

(b) How to formulate   and prove sampling theorems in the spirit of
\eqref{eqn:shannon}?

(c) Which sets in $G$ can take the role of $\ZZ \subseteq \RR$?  What
is ``uniform'' and ``nonuniform'' sampling in $G$?

These questions are of course interrelated and depend very much on
which notion of bandlimited functions on $G$ is chosen.

A natural attempt consists in replacing the Fourier transform on $\RR
$ by the operator-valued Fourier-Plancherel transform on $G$.  A
function is then said to be  be bandlimited, if  its
Fourier-Plancherel
transform is  supported in a given (quasi)-compact set of the dual
object $\widehat{G}$. This notion was pursued by Dooley for
motion groups, i.e.,  semidirect products of the form $G
= \RR^k \semdir K$, for a compact matrix group $K$~\cite{Do}.  He
derived uniqueness theorems for the resulting functions spaces  together
with reconstruction procedures.

In a different direction, Pesenson~\cite{Pe98} investigates sampling
problems on stratified nilpotent Lie groups $G$. In his context, the
Fourier transform on $\RR $ is replaced by the spectral decomposition
of the (left-invariant) sub-Laplacian on $G$. A function is then called
bandlimited, if it belongs to the range of a spectral projection of
the sub-Laplacian. Using a delicate optimization argument, he proved
uniqueness theorems for spaces on bandlimited functions.
 More generally, in a series of papers
\cite{Pe98,Pe01,Pes04}   Pesenson  proposed an
abstract notion of bandlimitedness associated to  any
(unbounded) self-adjoint operator on a Hilbert space and applied this
idea to a variety of situations.

Our contribution is threefold.

First, we develop the technique of
oscillations estimates for sampling in left-invariant closed subspaces
of $L^2(G)$ with a reproducing kernel. This technique has been
developed in the early 90's  to
treat nonuniform sampling of bandlimited functions and of
generalized wavelet transforms~\cite{FG92aa,Gr}. Recently oscillation estimates
have been rediscovered and adapted to sampling in shift-invariant
spaces~\cite{AF98},  in reproducing kernel Hilbert
spaces~\cite{FR05}, and sampling in Paley-Wiener spaces on
manifolds~\cite{FePe}.

Second, we solve the discretization problem for the continuous
wavelet transform with respect to discrete series
representations~\cite{ali-gazeau}. This is equivalent to a
sampling theorem for the generalized wavelet transform and amounts
to the construction of  frames contained in a single orbit of the
group action. Our contribution closes a technical gap in the
literature (where the discretization was shown under an additional
integrability condition~\cite{Gr}).

Third, we derive a sampling theorem for bandlimited functions on
stratified nilpotent Lie groups. The main theorem is a generalization
of a celebrated theorem of Beurling~\cite{beurling66} from $\RR ^n$ to
stratified groups.
 Though some of the proof ideas are similar in
spirit to Pesenson's~\cite{Pe98}, our approach yields a significant
technical simplification and several new insights. The stability of
the sampling procedure is made  explicit by formulating an analogue of
Shannon's theorem~\eqref{eqn:shannon} and~\eqref{eqn:shannon_norm}, and --- at least in
principle ---  the reconstruction of a bandlimited function from its
samples is stable and constructive.

\textbf{Notation.} In the following, $G$ denotes a
second-countable non-compact, locally compact group.   Then ${\rm
L}^2(G)$ is a separable Hilbert space.  We denote integration
against Haar measure by $\int_G \, \cdot \, dx$, and use $|A|$ to
denote the Haar measure of a Borel set $A \subset G$.
$\mathcal{U}_e$ denotes the neighborhood filter at unity.  The
left {\bf regular representation} $L$ of $G$ acts on ${\rm
L}^2(G)$ by $L_x f(y) = f(x^{-1}y)$. Given a function $f$ on $G$,
we let $f^*(x) = \overline{f(x^{-1})}$ be the usual involution.
 In the
following, we will use frequently that for $f,g \in {\rm L}^2(G)$,
the convolution  $f \ast g^*$ is a  well-defined  continuous function
vanishing at infinity, and that the convolution can be written as the
inner product  $f \ast g^*(x) = \langle f, L_x g \rangle$.

\vspace{3 mm}

\noindent \textbf{Acknowledgement.} Part of this work was carried out
when both authors were visiting the Erwin Schr\"odinger Institute in
Vienna. We would like to thank ESI for its great working conditions
and creative  research atmosphere.

\section{Sampling theorems and frames}

In this section we recall a few general properties of spaces with
sampling expansions and  observe  their close relationship to
reproducing kernel Hilbert spaces. This connection allows us to
view reconstruction formulas as a  discretization of the
reproducing kernel.

 The basic notion considered in this respect is that of a {\bf sampling
set}.
\begin{defn}
 Let ${\mathcal H} \subset {\rm L}^2(G)$ be a leftinvariant closed
 subspace  consisting of
 continuous functions. A discrete subset $\Gamma \subset G$ is called
 {\bf sampling set (for $\mathbf{\mathcal H}$)} if the restriction map $R_\Gamma : f \mapsto f|_\Gamma$ is
 a topological embedding ${\mathcal H} \to \ell^2(\Gamma)$, i.e., if there exist constants
 $0 < A \le B < \infty$ such that,  for all $f \in {\mathcal H}$,
 \begin{equation}
   \label{eq:c1}
 A \| f \|^2 \le  \sum_{\gamma \in \Gamma} |f(\gamma)|^2 \le B \| f \|^2~~.
 \end{equation}

When using the optimal constants in \eqref{eq:c1}, the quantity
$B/A$ is called {\bf tightness} of the sampling set.
\end{defn}

The sampling theorems in \cite{Do,Pe98} are  concerned with
establishing that $R_\Gamma$ is injective, whereas we are
interested in criteria to show that $R_\Gamma$ is a topological
embedding.

We first show that for left-invariant spaces
the {\em continuity} of the sampling process implies the existence
of  a reproducing kernel for ${\mathcal H}$ that is given by
convolution with a suitable ${\rm L}^2$-function.
\begin{thm} \label{eqn:rep_ker} Let ${\mathcal H} \subset {\rm L}^2(G) \cap C(G)$ be closed and
leftinvariant. Assume that there exists $\gamma \in G $ such that
the point evaluation  $f \mapsto f(\gamma)$ is a continuous linear
functional on ${\mathcal H}$. Then there exists a unique function
$p \in {\mathcal H}$ with the following properties: $p= p \ast
p^*$, and the orthogonal projection onto ${\mathcal H}$ is given
by the map $f \mapsto f \ast p^*$. Furthermore, $\cH $ consists of
continuous functions vanishing at infinity.
\end{thm}
\begin{prf}
 By the Riesz representation theorem there exists $p_\gamma \in \cH$ such
that $f(\gamma) = \langle f, p_\gamma \rangle$ holds for all $f
\in {\mathcal H}$. We let $p = L_{\gamma^{-1}} p_\gamma$. Since
${\mathcal H}$ is leftinvariant, $L_{\gamma x^{-1}} f \in
{\mathcal H}$, and we obtain for all $x \in G$ that
\begin{eqnarray*}
 f(x) & = & f(x\gamma^{-1} \gamma) = (L_{\gamma x^{-1}} f) (\gamma) =
\langle L_{\gamma x^{-1}} f, p_\gamma \rangle = \langle f, L_x
(L_\gamma^{-1} p_\gamma) \rangle \\ & = & \langle f, L_x p \rangle
= f \ast p^*(x)~~.
\end{eqnarray*}
Consequently, if $f\in \cH$, then $f = f \ast p^*$, and if $f\in \cH
^\perp $, then $f \ast p^* = 0$.
Thus  the orthogonal projection onto ${\mathcal H}$ is indeed
given by right convolution with $p^*$ and $\cH $ consists of
continuous functions vanishing at infinity.

Inserting $f=p$, we obtain $p = p \ast p^*$, and thus $p^* = (p
\ast p^*)^* 
 = p \ast p^* = p$. For uniqueness,
assume that $q \in {\mathcal H}$ with $q = q \ast q^*$, and such
that $f \mapsto f \ast q^*$ is the projection onto ${\mathcal H}$.
Observe that then $q = q^*$ and $p = p^*$, and therefore
\[ q = q^* = (q \ast p)^* = p^* \ast q^* = p \ast q = p ~~.\]
\end{prf}


 Recall that a {\bf frame} in a Hilbert space ${\mathcal
H}$ is a family $(\eta_i)_{i \in I}$ satisfying for all $f \in
{\mathcal H}$ the inequalities
\[
 A \| f \|^2 \le \sum_{i \in I} |\langle f, \eta_i \rangle |^2 \le B \| f \|^2~~,
\]
for constants $0 < A \le B < \infty$.

It has been observed on several occasions that sampling theorems in reproducing
kernel Hilbert spaces  are closely related to frames. The
following result formulates this for the group case.
\begin{cor} \label{cor:sample_frame}
 Let ${\mathcal H} \subset {\rm L}^2(G)$ a leftinvariant
 closed subspace consisting of continuous functions, and such that point evaluations
are continuous on ${\mathcal H}$. Let $p$ be the associated
reproducing kernel of ${\mathcal H}$. Then the following are
equivalent.
 \begin{itemize}
 \item[(a)] $\Gamma$ is a sampling set.
 \item[(b)] The family $(L_\gamma p)_{\gamma \in \Gamma}$ is a frame of ${\mathcal
H}$.
 \end{itemize}
\end{cor}
\begin{prf}
 The statement follows from the fact that
 sampled values and frame coefficients coincide.
\end{prf}

\begin{rem}
The previous observations allow us to apply well-known results from
frame theory to  obtain sampling expansions. In particular, the
reconstruction of a function  from its samples can be achieved by the
frame algorithm~\cite{DS52,Chr03}. Let
\begin{equation}
  \label{eq:c6}
Sf =  \sum _{\gamma \in \Gamma} \langle f , L_\gamma p \rangle L_\gamma
p =  \sum _{\gamma \in \Gamma}  f (\gamma )  L_\gamma p
\end{equation}
be the frame operator associated to the family $\{ L_\gamma p : \gamma
\in \Gamma \}$. Then $\langle Sf, f \rangle = \sum _{\gamma \in \Gamma
  } |f(\gamma )|^2 \geq A \|f\|^2$ for $f\in \cH $, and so $S$ is
  invertible on $\cH $. Set $\tilde{e_\gamma } = S\inv (L_\gamma p)$,
  the so-called \emph{dual frame}, then
  \begin{equation}
    \label{eq:c7}
    f = S\inv S f = S\inv \big( \sum _{\gamma \in \Gamma } f(\gamma)
    L_\gamma p \big) = \sum _{\gamma \in \Gamma } f(\gamma)
    \tilde{e_\gamma } \, ,
  \end{equation}
and so we have a sampling expansion that is in complete analogy to the
cardinal series~\eqref{eqn:shannon}.
 The above series (\ref{eq:c7}) converges unconditionally in
 ${\rm L}^2 (G)$ by frame theory,
but it  also converges  {\em uniformly}. Since  $f (x) =
(f \ast p^*)(x) = \langle f, L_x p \rangle$ and thus $|f(x)| \le \| f
\| \|p \|$ uniformly in $x\in G$, we obtain the  uniform
convergence of \eqref{eq:c7} from the $L^2$-convergence.

The reconstruction formula can be  simplified in the following two
cases.

If the sampling set is tight, i.e. $A=B$ in \eqref{eq:c1}, then
the frame operator $S$ is $A \cdot {\rm Id}_{\cH }$, and so
$\tilde{e_\gamma} = S\inv (L_\gamma p)  = \frac{1}{A} L_\gamma p$.
The reconstruction $f= \frac{1}{A}\sum _{\gamma \in \Gamma }
f(\gamma ) L_\gamma p$ is then the exact analog of
\eqref{eqn:shannon}.

Furthermore, if the sampling set  $\Gamma$ is a subgroup, then $S$
commutes with the translations $L_\gamma , \gamma \in \Gamma $, and so
$S\inv (L_\gamma p) = L_\gamma S\inv p = L_\gamma \tilde{e}$ for a
single function $e\in \cH $. As a result we obtain the following
analog of Shannon's sampling theorem:
\begin{equation} \label{eqn:rec_sampl_lattice}
 f(x) = \sum_{\gamma \in \Gamma} f(\gamma) (L_\gamma \tilde{p})(x)~~.
\end{equation}
This reconstruction formula is in nice correspondence  to the continuous
 formula
\begin{equation} \label{eqn:rec_cont}
f(x) = \langle f, L_x p \rangle = \int_G f(y) (L_y p^*)(x) dy~~.
\end{equation}
Now (\ref{eqn:rec_sampl_lattice}) is easily recognized as a
Riemann sum analog of (\ref{eqn:rec_cont}). Hence one expects that
$\tilde{p} \approx c_{\Gamma} p$, as the tightness of the frame
approaches one.
\end{rem}

\begin{rem} In the light of Theorem \ref{eqn:rep_ker}, Dooley's
  definition of bandlimited functions turns out to be  too large for
  the existence of Shannon-type sampling theorems. Using the
Plancherel formula for the motion group, one can show that Dooley's
spaces  have an ${\rm L}^2-$reproducing kernel only  if the group
$G$ is a {\em finite} extension of a vector group. See
\cite{Fu_LN} for more details.
\end{rem}

\section{Sampling and oscillation estimates}
The sampling theorems obtained in this paper are derived using
${\rm L}^2$-estimates of the oscillatory behavior of the functions
under consideration.  Following~\cite{Gr}, we define the
modulus of continuity, the so-called \emph{oscillation}, of a function
$f$ on $G$.
\begin{defn}
 For a function $f$ on $G$ and a set $U \subset G$, we define
 \[
{\rm osc}_U(f) (x) = \sup_{y \in U} |f(x)-f(xy^{-1})|~~.
 \]
\end{defn}
Clearly, $U \subset U'$ implies ${\rm osc}_U(f) \le {\rm
osc}_{U'}(f)$ pointwise. Moreover, if $U$ is relatively compact
and $f$ is continuous, then ${\rm osc}_U(f)$ is
well-defined everywhere and continuous: Indeed, for $x,z \in G$,
\begin{eqnarray*}
|{\rm osc}_U(f)(x) - {\rm osc}_U(f)(z)| & = & \left| \sup_{y \in
U} |f(x)
- f(xy^{-1}) | - \sup_{y \in U} |f(z) - f(zy^{-1}) |  \right| \\
& \le & \sup_{y \in U} \left| \raisebox{0cm}[0.2cm][0.2cm]{} |f(x)
- f(xy^{-1})| - |f(z) - f(zy^{-1})| \right| \\ & \le & |f(x)-f(z)|
+ \sup_{y \in U} |f(xy^{-1})-f(zy^{-1})|~~.
\end{eqnarray*}
If  $z \to x$, then the  first term tends to zero by continuity of
$f$.  The second term goes to zero, because $f$ is (left)
uniformly  continuous on the relatively compact set $xU$.


We next consider conditions on the sampling set. These
requirements are quite intuitive and generalize certain density
concepts from $\RR ^n $ to $G$.

Conditions of this type are often encountered in connection with
discretization,  see e.g. \cite{beurling66,FeiGr1}.
\begin{defn} Let $U,W \subset G$ be Borel, and $\Gamma \subset G$ countable.\\
(a) $\Gamma$ is called {\bf $\mathbf{U}$-dense} if $\Gamma U =
\cup_{\gamma \in \Gamma} \gamma U = G$. \\ (b) $\Gamma$ is called
{\bf $\mathbf{W}$-separated} if $|\gamma W \cap \gamma' W| = 0$
for distinct $\gamma,\gamma' \in \Gamma$. $\Gamma$ is called {\bf
separated} if it is $W$-separated for some $W \in {\mathcal U}_e$.
\\ (c) $\Gamma$ is a {\bf quasi-lattice} if there exists a
relatively compact Borel set $C$, such that $\Gamma$ is both
$C$-separated and $C$-dense. Such a set $C$ is called a {\bf
complement of $\mathbf{\Gamma}$}.
\end{defn}

The main example of a quasi-lattice is a cocompact, discrete
subgroup $\Gamma < G$ (often called a lattice in $G$).  Here any
relatively compact fundamental domain of $\Gamma $, i.e. a
relatively compact set of representatives mod $\Gamma$, can be
chosen as  a complement $C$. However, the concept of a
quasi-lattice is strictly weaker than that of a cocompact discrete
subgroup, and quasi-lattices may exist even in groups that do not
admit a lattice.

We next collect some technical lemmas in connection with separated
and dense sets.
\begin{lemma} \label{lem:ex_sample}
 For every $U \in {\mathcal U}_e$ there exists a separated
$U$-dense set.
\end{lemma}
\begin{prf}
 Let $V \in {\mathcal U}_e$ satisfying $V V^{-1} \subset U$.
Choose  a $V$-separated set $\Gamma$ that is maximal with respect to
inclusion. Then  by
maximality, for every $g\in G$  there exists $\gamma \in \Gamma$ such that
$g V \cap \gamma V \not= 0$. But this implies $g \in \gamma
VV^{-1} \subset \gamma U$, and so $\Gamma $ is $U$-dense.
\end{prf}

The following lemma provides a  substitute for a uniform partition of
$G$. 

\begin{lemma} \label{lem:well-spread}
 Let $U,W \in \cU _e$ with $W\subset U$, and  $\Gamma \subset G$ be a
 $W$-separated and $U$-dense set.  Then there exist
relatively compact Borel sets $(V_\gamma)_{\gamma \in \Gamma}$
such that  $W \subset V_\gamma \subset U$ and $(\gamma
V_\gamma)_{\gamma \in \Gamma}$ is a partition of $G$, i.e., $G =
\bigcup_{\gamma \in \Gamma}^\bullet \gamma V_\gamma$ as a {\em
disjoint} union.
\end{lemma}

\begin{prf}
Since $G$ is $\sigma$-compact, the cardinality of disjoint
translates of $W\in \mathcal{U}_e$ can be at most countable. Since
$\Gamma $ is $W$-separated and $U$-dense, $\Gamma $ must be
countable.  We may therefore enumerate $\Gamma$ using a bijection
$\NN \to \Gamma$ and write  $\Gamma = \{\gamma_k : k \in \NN \}$.
Now let $A = G \setminus \bigcup_{\gamma \in  \Gamma} \gamma W$,
which is a Borel set. We define recursively
\begin{equation}
  \label{eq:c2}
 V_{\gamma_k} = W \cup \gamma_{k}^{-1} \left( (A \cap \gamma_k U)
\setminus \bigcup_{i=1}^{k-1} \gamma_i V_{\gamma_i} \right)~~.
\end{equation}
Observe that the union is {\em disjoint}, since $\gamma_k W
\subset G \setminus A$, whereas \[ \gamma_k \gamma_{k}^{-1} \left(
(A \cap \gamma_k U) \setminus \bigcup_{i=1}^{n-1} \gamma_i
V_{\gamma_i} \right) \subset A~~.
\]
 We claim that the $V_\gamma$'s
have the desired properties. Clearly, $W \subset V_\gamma \subset U$
and the $V_{\gamma _j}$'s are  measurable by construction. Now suppose
that $g \in
\gamma_j V_{\gamma_j} \cap \gamma_k V_{\gamma_k}$ for $j<k$. On the
one hand,   if $g \in \gamma_k W$, then $g \in G \setminus A$, whereas
$g \in \gamma_j V_j$ implies  $g \in \gamma_j W$. This
contradicts the $W$-separatedness.  On the other hand, if $g \in A$,
then
\[ g \in (A \cap \gamma_k U) \setminus \bigcup_{i=1}^{n-1} \gamma_i
V_{\gamma_i} ~~,\] which contradicts $g \in \gamma_j V_j$, as $j
<k$.

Finally, let $g \in G$ be arbitrary. If $g \in \Gamma W$, then $W
\subset V_\gamma$ shows $g \in \gamma V_\gamma$ for a suitable
$\gamma$. In the other case let $k \in \NN$ be minimal with $k \in
\gamma_k U$. Then $g \in \gamma_k V_{\gamma_k}$, since
\[\gamma_k V_k \supset (A \cap \gamma_k U) \setminus
\bigcup_{i=1}^{k-1} \gamma_i V_{\gamma_i} \supset (A \cap \gamma_k
U) \setminus \bigcup_{i=1}^{k-1} \gamma_i U ~~.\] By assumption on
$k$, $g$ is contained in $\gamma_k U$, but not in $\gamma_i U$,
for $i<k$, which proves the statement.
\end{prf}

The next theorem shows how  to derive sampling theorems from
oscillation estimates.
\begin{thm} \label{thm:osc_sampl}
 Let $U,W \in \cU _e$ and $W\subset U$ and  $\Gamma \subset G$ be  $U$-dense and
$W$-separated set.   Let
${\mathcal H} \subset {\rm L}^2(G)$ be a closed, leftinvariant
subspace of ${\rm L}^2(G)$  consisting of continuous functions and assume
that there exists $\epsilon, 0 < \epsilon <1$, such that  $\| {\rm osc}_U (f) \|_2 \le \epsilon \| f
\|_2$ for all $f \in {\mathcal H}$. Then $\Gamma$ is a sampling set for ${\mathcal H}$. More
precisely, we have the estimate
\begin{equation}
\label{eqn:sampl_equiv} \frac{1}{|U|^2} (1-\epsilon)^2 \| f \|^2_2
\le \| R_\Gamma f \|^2_2 \le \frac{1}{|W|^2} (1+\epsilon)^2 \| f
\|^2_2~~, \qquad \forall f \in {\mathcal H}
\end{equation}
\end{thm}

\begin{prf} Let $(V_\gamma)_{\gamma \in \Gamma}$ be the family of Borel sets
 asserted  by Lemma~\ref{lem:well-spread}. We
introduce the auxiliary operator $Q: \ell^2(\Gamma) \to {\rm
L}^2(G)$  defined by
\[
 Q((c_\gamma)_{\gamma \in \Gamma}) = \sum_{\gamma \in \Gamma} c_\gamma L_\gamma \chi_{V_\gamma}~~.
\]
Since the sets $\gamma V_\gamma$ are  pairwise disjoint and  $|W| \le
|V_\gamma| \le |U|$, $Q$ is
a well-defined bounded operator with operator norm $\| Q \|_\infty
= \sup_{\gamma \in \Gamma} |V_\gamma| \le |U|$. More importantly,  $Q$ has a
bounded inverse on its range  with operator norm  $\| Q^{-1} \|_\infty =
\frac{1}{\inf_{\gamma \in \Gamma} |V_\gamma|} \le \frac{1}{|W|}$.

Since $Q$ is an interpolation of the sequence $(c_\gamma )$ by a step
function, it is plausible that  $Q$ approximates  well  the inverse of
$R_\Gamma$. The following estimate makes this precise by introducing
${\rm osc}_U$ to the argument. Since  $V_\gamma \subset U$ and the
$\gamma V_\gamma$'s are disjoint, we may estimate, for all $f \in
{\mathcal H}$,
\begin{eqnarray*}
\| f - Q R_\Gamma f \|^2_2 & = & \sum_{\gamma \in \Gamma}
\int_{\gamma V_\gamma} |f(x) - f(\gamma)|^2 dx
\\ & \le & \sum_{\gamma \in \Gamma} \int_{\gamma V_\gamma} |{\rm osc}_U (f) (x)|^2 dx \\
& = & \| {\rm osc}_U(f)\|^2_2 \\
& \le & \epsilon^2 \| f \|^2_2~~.
\end{eqnarray*}
Consequently, we obtain the upper bound of the sampling inequality
\eqref{eq:c1} for  $f \in {\mathcal H}$
\begin{eqnarray*}
 \| R_\Gamma f \|_2 & = &  \| Q^{-1} Q R_\Gamma f \|_2 \\
& \le & \| Q^{-1} \|_\infty \| Q R_\Gamma f \|_2 \\
& \le & \| Q^{-1} \|_\infty (\| f \|_2 + \| f - Q R_\Gamma f \|_2) \\
& \le & \| Q^{-1} \|_\infty (1+\epsilon) \| f \|_2 \\
& \le & \frac{1}{|W|} (1+\epsilon) \| f \|_2~~.
\end{eqnarray*}
The decisive lower bound follows similarly by
\begin{eqnarray*}
 \| R_\Gamma f \|_2 & \ge &  \| Q \|_\infty^{-1} \| Q R_\Gamma f \|_2 \\
& \ge & \| Q \|_\infty^{-1} (\| f \|_2 - \| f - Q R_\Gamma f \|_2) \\
& \ge & \| Q \|_\infty^{-1} (1-\epsilon) \| f \|_2 \\
& \ge & \frac{1}{|U|} (1-\epsilon) \| f \|_2~~.
\end{eqnarray*}
Thus $\Gamma $ is a sampling set for $\cH $.
\end{prf}

\begin{rem}
The theorem allows to estimate the tightness of the sampling
estimate from above by
\[
\frac{|U|^2}{|W|^2} \cdot \frac{(1+\epsilon)^2}{(1-\epsilon)^2}~~.
\]
 The first quotient $|U|^2/|W|^2$  is a measure
for the  {\bf uniformity} of the sampling set $\Gamma$. For
quasi-lattices it  is $|U|^2/|W|^2 = 1$ and we may therefore  call
this case  {\em uniform sampling}. In this case the tightness
estimate in Theorem \ref{thm:osc_sampl} depends only  on the
oscillation.

The second quotient $\frac{(1+\epsilon)^2}{(1-\epsilon)^2}$
depends on the {\bf density} of $\Gamma$. Since  ${\rm osc}_U (f)
\to 0$ as $U \to \{e\}$, high density (small $U$) improves the
tightness of the sampling procedure.
\end{rem}

Theorem \ref{thm:osc_sampl} and Lemma \ref{lem:ex_sample} yield
the following existence result.\begin{cor} \label{cor:ex_ss}
 Let ${\mathcal H}$ be a leftinvariant closed space consisting of
continuous functions that satisfy $\| {\rm osc}_U (f)  \|_2 \le
\epsilon \| f \|_2$ for $\epsilon < 1$ and a suitable $U \in
{\mathcal U}_e$. Then there exists a sampling set for ${\mathcal
H}$.
\end{cor}

A simple but useful trick for the derivation of oscillation
estimates is the following observation:
\begin{equation} \label{eqn:osc_conv} \begin{split}
 {\rm osc}_U(f \ast g)(x)  & = \sup_{z \in U} |\int_G f(y) (g(y^{-1} x) -g(y^{-1}xz^{-1})
dy|\\ & \le  \int_G |f(y)| \sup_{z \in U} |g(y^{-1} x) -
g(y^{-1}xz^{-1})| dy\\ & \le  (|f| \ast {\rm osc}_U(g))
 (x)~~. \end{split}
\end{equation}
This can be combined with Theorem \ref{thm:osc_sampl} to establish
sampling theorems for a certain class of leftinvariant spaces.
Similar arguments were  employed  in \cite{FeiGr2,Gr}.

\begin{thm} \label{thm:sampl_cv1}
Let ${\mathcal H}$ be a closed, leftinvariant subspace of $L^2(G)$. Assume
that there exists a continuous $g \in {\rm L}^1(G)$ satisfying $\|
{\rm osc}_W g \|_1 < \infty$ for some $W \in {\mathcal U}_e$, as
well as $f = f \ast g$ for all $f \in {\mathcal H}$. Then there
exists $U \in {\mathcal U}_e$ such that every separated $U$-dense
set is a sampling set.
\end{thm}

\begin{prf}
The continuity of $g$ yields that ${\rm osc}_U (g) \to 0$
pointwise, as $U$ runs through a basis of ${\mathcal U}_e$. Since
$\| {\rm osc}_W g \|_1 < \infty$, and ${\mathcal U}_e$ has a
countable basis, the dominated convergence theorem applies to
yield $ {\rm osc}_U g \to 0$ in the ${\rm L}^1$-norm. Pick $U$
with $\| {\rm osc}_U g \|_1 < 1$. Then inserting
(\ref{eqn:osc_conv}) into Theorem \ref{thm:osc_sampl} yields the
desired statement.
\end{prf}

\section{Sampling Theorems and Discrete Series Representations}

In this section we consider a particular case of left-invariant closed
subspaces of $L^2(G)$ that arise in the context of the representation
theory of $G$. The  reproducing kernel Hilbert  spaces
 described in Theorem \ref{eqn:rep_ker} occur naturally as the range
 of a  (generalized) continuous wavelet transforms. These
 are obtained by the following procedure.  Given a unitary
 representation $(\pi, {\mathcal H}_\pi)$ and a vector $\eta \in
 {\mathcal H}_\pi$, we define the (generalized) \emph{wavelet transform}
 $V_\eta$ from $\mathcal{H}_\pi $ to $L^\infty (G)$ by
$$ V_\eta \varphi (x) =  \langle \varphi, \pi(x) \eta \rangle \,
\qquad x\in G \, .
$$
This operator maps vectors $\phi \in \cH _\pi $ onto
representation coefficients of $G$.
 We call $\eta$
 {\bf admissible} whenever $V_\eta$ is an isometry into $L^2(G)$. The
 properties  of the regular representation of
 $G$ lead to the following conclusions: (a)  the space  ${\mathcal H}=
 V_\eta ({\mathcal H}_\pi)$ is   a
 closed, leftinvariant subspace of $L^2(G)$, (b) $V_\eta $ intertwines
 $(\pi ,\cH )$ and the regular representation $L_x$ restricted to  $V_\eta (L^2(G))$,
 and (c)
 the projection from $L^2(G)$ onto ${\mathcal H}$ is given by right
 convolution with  $S = V_\eta
 \eta$.  See e.g. \cite{Fu_LN} and the references therein  for
 details. Now the construction  of a sampling set
 for ${\mathcal H}$ is equivalent to  the problem of finding
 $\Gamma \subset G$ such that $\pi(\Gamma) \eta$ is a frame of
 ${\mathcal H}_\pi$. This question is referred to as the {\em
   discretization problem}  for the continuous wavelet transform, and
 has attracted considerable attention~\cite{ali-gazeau,Chr03}.

A special class of representations for which the construction of
frames has been investigated extensively are the  so-called {\bf
discrete series representations}, i.e., irreducible
square-integrable representations of $G$. These always possess
admissible vectors. The papers
\cite{FeiGr1,Da,Gr,BeTa,Fu,AnCaDeLe}  are a small, but
non-exhaustive list  of papers where the discretization problem of
discrete series representations has been studied.

To our knowledge, the construction of frames from discrete series
representations has been proven rigorously only under the additional
assumption that the representation be \emph{integrable}.
 The following result shows that {\em all}
discrete series  representations can be discretized and yield a
construction of frames. Despite
the widespread interest in this question, the result seems to be
new (although it has been mentioned in \cite{Gr}).

\begin{thm} \label{thm:disc_ser}
 Let $\pi$ be a discrete series representation of $G$. There exists
 a  vector $\eta \in {\mathcal H}_\pi$ and $U \in {\mathcal U}_e$
 (depending on $\eta $ and $\pi $) such that
 $(\pi(\gamma) \eta)_{\gamma \in \Gamma}$ is a frame for  ${\mathcal H}_\pi$, whenever
 $\Gamma \subset G$ is   separated  and  $U$-dense.
\end{thm}

\begin{prf}
 Let $\psi$ be an arbitrary admissible vector,
 ${\mathcal H}_0 = V_\psi({\mathcal H}_\pi)\subseteq L^2(G)$ and $P$
 be the orthogonal projection of $L^2(G)$ onto $\cH _0$. Choose
 $h \in   C_c(G)$ that projects onto a nonzero element of ${\mathcal
   H}_0$. This means that $Ph = V_\psi \eta $ for some non-zero $\eta
 \in \cH _\pi $.  Since $V_\psi $ is an isometry and
 intertwines with the regular representation of $G$, we obtain 
 \begin{eqnarray*}
 V_\eta \varphi (x) & = &  \langle \varphi, \pi (x) \eta \rangle  \\
& = & \langle  V_\psi \varphi, L_x V_\psi \eta \rangle \\
&=& \langle P  V_\psi \varphi, P (L_x  h)  \rangle \\
& = & \langle V_\psi \varphi,  L_x h \rangle \\
 & = & \left( V_\psi \varphi  \ast h^* \right)(x)\, . ~~
 \end{eqnarray*}
Therefore $V_\eta (\pi (z)\phi) = V_\psi (\pi (z) \phi ) \ast h^* =
L_z (V_\psi \phi \ast h^*) = L_z V_\eta \phi$, and  $V_\eta$ is a
(nonzero) bounded intertwining operator from
 ${\mathcal H}_\pi$ into ${\rm L}^2(G)$. 
 Since $\pi$ is  irreducible and $\eta$ is nonzero, Schur's lemma
 implies that $V_\eta$ is a  scalar multiple of an isometry and that
 $V_\eta $  maps  ${\mathcal H}_\pi$ onto a closed leftinvariant
 subspace ${\mathcal H} \subset {\rm L}^2(G)$.

 We can now  employ
 (\ref{eqn:osc_conv}) to establish
 \[ \| {\rm osc}_U V_\eta \varphi
 \|_2 \le \| V_\psi \varphi \|_2 \| {\rm osc}_U h^* \|_1 = c_{\eta,\psi} \|
 V_\eta \varphi \|_2 \| {\rm osc}_U h^* \|_1 ~~,\] for a suitable
 positive constant $c_{\eta,\psi}$.
Since $h\in C_c(G)$, the oscillation  ${\rm osc}_U h^*$ is a
bounded and  compactly supported for all relatively compact sets
$U$.  In particular,  ${\rm osc}_U h^*$ is integrable, and by
choosing  $U$ small enough, we obtain $c_{\eta,\psi} \| {\rm
osc}_U h^* \|_1< 1$.  This is exactly, what is needed to apply
Theorem \ref{thm:osc_sampl}. We conclude that every separated and
$U$-dense set $\Gamma \subset G$ is a sampling set for $ \cH $,
and $(\pi(\gamma) \eta)_{\gamma \in \Gamma}$ is a frame for $\cH
_\pi $.
\end{prf}

\begin{rem}
  It is easy to see that the set of $\eta $ for which the above
  argument works, is a dense subspace of $\cH _\pi$. If in addition,
  $\pi $ is integrable, then a different argument yields the existence
  of ``frame vectors'' $\eta $~\cite{Gr}.

  For nilpotent connected Lie groups, all  discrete series representations
  are in fact integrable; this follows by \cite[Theorem 4.5.11]{CoGr}.
  However, the semisimple Lie group $SL(2,\RR)$ provides
  an example of a discrete series representation that is not
  integrable.
\end{rem}



\section{Paley-Wiener spaces on stratified Lie groups}

As a second application we use the oscillation method of Section~3
to derive sampling theorems on stratified Lie groups.
In the setting of stratified Lie groups, the required oscillation
estimates can be formulated and proved with particular ease. The
intuition behind this approach is that the control over the
sub-Laplacian entails control over all finite order left-invariant
differential operators. Hence the oscillation can be controlled.

We first recall some basic facts about stratified Lie groups and then
define Paley-Wiener spaces on
 such groups. We refer to \cite{CoGr,FoSt,Fo75} for more
details. We assume that $G$ is a simply connected, connected
nilpotent Lie group with Lie algebra $\mathfrak{g}$ of dimension
$n$. The Lie algebra is assumed to be {\bf stratified}, which
means that $\mathfrak{g}$ is the direct sum of subspaces
$V_1,\ldots,V_m$ satisfying $[V_1,V_j] = V_{j+1}$ for $1 \le j \le
m$, where we use $V_{m+1} = \{ 0 \}$. The {\bf homogeneous
dimension} of $G$ is given by $Q = \sum_{j =1}^m j \cdot ({\rm
dim}V_j)$. We will use results from \cite{Fo75}, and therefore we
adopt the assumption $Q>2$ made in that paper, noting that it only
excludes the groups $\RR$ and $\RR^2$, for which the sampling
theorems derived below are known anyway. For the following, we fix
a basis $X_1,\ldots,X_n$ of $\mathfrak{g}$ that is composed of
bases of the $V_j$, i.e., $X_k \in V_j$  for $\sum_{i=1}^{j-1}
{\rm dim}(V_i) < k \le \sum_{i=1}^{j} {\rm dim}(V_i)$.

If $\mathfrak{g}$ is  stratified, it possesses a
one-parameter group of Lie algebra automorphisms defined as
\[
\delta_t (\sum_{i=1}^m v_i) = \sum_{i=1}^m t^i v_i ~~,v_i \in
V_i~~.
\]
We also fix a {\bf homogeneous norm}, which is a mapping $| \cdot | : \mathfrak{g} \to \RR^+$,
fulfilling $|\delta_t(X)| = t |X|$ and $|-X| = |X|$. Confer
\cite{FoSt} for existence.


Since for simply connected, connected nilpotent Lie groups the
exponential map ${\rm exp}: \mathfrak{g} \to G$ is a polynomial
diffeomorphism with polynomial inverse, the dilations $\delta _t$
on $\mathfrak{g}$ yield a one-parameter group of  automorphisms
$\delta _t$ of $G$, and a homogeneous norm $|.|$ on $G$. (As is
costumary,  we use the same notation on $\mathfrak{g}$ and on
$G$.)   The homogeneous norm fulfills $|\delta_t(x)| = t |x|$,
$|x^{-1}|=|x|$, and  the triangle inequality
\begin{equation} \label{eqn:triang} |xy| \le C_\triangle(|x| + |y|)
  \quad \quad x,y \in G
~~,\end{equation} for some  constant
$C_\triangle>0$~\cite[Proposition 1.6]{FoSt}.  The Haar measure is
changed by the dilations $\delta _t$ according to the formula
$|\delta_t(A)|= t^Q |A|, A\subseteq G$. Consequently, on ${\rm
L}^2(G)$ we have
\begin{equation} \label{eqn:dil_L2}
\| f \circ \delta_t \|_2 = t^{-Q/2} \| f \|_2~~.
\end{equation}

Next we consider differential operators on $G$. $C^{k}(G)$ denotes
the space of $k$ times continuously differentiable functions on
$G$. We identify $\mathfrak{g}$ in the usual manner with the space
of leftinvariant differential operators of order one  acting on
$C^{\infty}(G)$. We use the multiindex notation $X^\alpha$, for
$\alpha \in \NN_0^{n}$, to denote the monomial differential
operator $X_1^{\alpha_1}\cdots X_n^{\alpha_n}$ of order $|\alpha|
= \sum_{i=1}^m |\alpha_i|$, where the $X_i$ are the elements of
the basis fixed above. Since  $\delta_t$ acts on $V_j$ by
multiplication with $t^j$,  all  $X_i$ are homogeneous, and
\begin{equation}
  \label{eq:c8}
 X_i (f  \circ \delta_t)  = t^j (X_i f) \circ \delta_t \, \quad
 \mathrm{for} \,\,\, X_i \in V_j, \forall f \in C^\infty(G) \, .
\end{equation}
 Consequently, all monomial
differential operators $X^\alpha$
inherit a similar homogeneity property, \begin{equation}
\label{eqn:homogen} X^\alpha ( f \circ \delta _t)  = t^{d(\alpha)}
(X^\alpha f) \circ \delta_t~~,
\end{equation} where $d(\alpha)$ is a suitable integer $\ge |\alpha|$.

For the analysis on stratified Lie groups the so-called  {\bf
  sub-Laplacian} $\mathcal{L}$ plays a distinguished role. If  $\ell
= {\rm dim}(X_1)$, then $X_1,\ldots,X_\ell$ is a basis of $V_1$,
and the sub-Laplacian is defined as $\mathcal{L}=-\sum_{i=1}^\ell
X_i^2$. It is well-known that $\mathcal{L}$ extends uniquely from
$C_c^\infty(G)$ to a selfadjoint positive definite operator on
${\rm L}^2(G)$. By the spectral theorem, we can associate to
$\mathcal{L}$ a projection-valued measure (or  spectral measure),
which we denote by $\Pi_\mathcal{L}$. Following
Pesenson~\cite{Pe98},   the {\bf
  Paley-Wiener space}
$E_\omega (G)$, for $\omega>0$, is defined by $E_\omega =
\Pi_\mathcal{L}([0,\omega]) ({\rm L}^2(G))$. Since $\mathcal{L}$ is left-invariant, the
projections $\Pi _\mathcal{L} ([0,\omega ])$ are also left-invariant, and
therefore 
$E_\omega $ is a closed, leftinvariant subspace of ${\rm L}^2(G)$.

\begin{Ex}
  If $G = \RR ^n$, then the sub-Laplacian is simply the Laplace
  operator $-\Delta $. The Fourier transform yields the spectral
  representation $(-\Delta f)\, \widehat \, (\xi ) = |\xi |^2 \hat{f}
  (\xi )$, therefore $(\Pi _\Delta ([0,\omega ]) f)\, \widehat{} \,
  (\xi ) = \chi _{\{\xi : |\xi |^2 \leq \omega \}}|\xi |^2 \hat{f}
  (\xi )$, and so the Paley-Wiener space $E_\omega $ is identical to
  the space of ``bandlimited functions'' $\{ f\in L^2(\RR ^n):
  \mathrm{supp}\, \hat{f} \subseteq B_{\sqrt{\omega }}\}$.

This example illustrates that $E_\omega (G)$ is a reasonable
generalization of the
usual notion of Paley-Wiener space on $\RR ^n$.
\end{Ex}

The following  theorem summarizes the main property of $E_\omega$
that is required for the derivation of sampling theorems.

\begin{thm}\label{thm:bernstein} For all $\omega>0$, $E_\omega \subset C^\infty(G)$.
For every $\alpha \in \NN_0^n$, the differential operator
$X^{\alpha}: E_\omega \to {\rm L}^2(G)$ is bounded.
\end{thm}
\begin{prf}
 The domains of the powers of the sub-Laplacian define a scale of Sobolev spaces
$S_s^2$ ($s>0$) \cite{Fo75}. The Bernstein inequality
\cite[Theorem 1]{Pe98} states for all $f \in E_\omega$ and all
$k\in \NN$ that $\| \mathcal{L}^k f\| \le \omega^k \| f \|$, in
particular $E_\omega \subset S_{2k}^2$. This holds for all $k\ge
0$, hence $E_\omega \subset C^\infty(G)$.
 Now by Corollary 4.13
of \cite{Fo75}   the Sobolev norm,  for fixed $k \in \NN$,
\[ \| f \|_{2,2k} = \| f \|_2 + \| {\mathcal L}^{k} f
\|_2 \] is equivalent to the norm
\[ f \mapsto \sum_{|\alpha| \le 2k} \| X^\alpha f \|_2 ~~.\]
Consequently,  the Bernstein inequality implies for all $f \in E_\omega$
\[ \sum_{|\alpha| \le 2k} \| X^\alpha f \|_2 \le C \omega^k \| f
\|_2 ~~,\] with a constant $C$ independent of $f$.
\end{prf}

Our goal is to establish a sampling theorem for the Paley-Wiener space
$E_\omega (G)$. In view of Theorem~\ref{thm:bernstein}
the basic strategy is to 
 derive a uniform ${\rm L}^2$-estimate for the  oscillation of $f\in
 E_\omega $. Applying suitable dilations $\delta _t $ to a $U$-dense and
 $W$-separated set   $\Gamma \subseteq G$, we can produce
a set $\delta _t(\Gamma )$ of any required density, while
preserving the uniformity. It is then plausible that the tightness
of a sampling estimate improves with increasing  density.


Therefore we try to derive  estimates of ${\rm osc} (U)$ as a
function of the diameter of $U$. As a tool we will use the  mean
value theorem, which is the simplest version  of Taylor's formula
for $G$, and  Sobolev-type estimates. We  cite  the mean value
theorem, which is a left-invariant version of \cite[1.33]{FoSt}.

\begin{lemma} \label{thm:mean_value}
 Let $G$ be stratified. There exist constants $C>0$ and $b\geq 1$ such
 that for all
 $f \in C^1(G)$ and all $x,y \in G$,
 \[
  |f(xy)-f(x)| \le C |y| \sup_{|z| \le b |y|,1\le j \le n} |X_j f(xz)|~~.
 \]
\end{lemma}

Next we state  a Sobolev-type estimate for  the comparison of a
local uniform norm and the  ${\rm L}^2$-norm.   Given any function
$f$ on $G$, $U \subset G$, we write  $\| f \|_{p,U} = \| f \cdot
\chi_U \|_p$ for the local ${\rm L}^p$-norms.

\begin{lemma}[\cite{CoGr}, Lemma A.1.5]
 For each compact set  $K \subset G$ there is a constant $C_K$ such that
 \[
  \| f \|_{\infty, K} \le C_K \sum_{|\alpha| \le n} \| X^\alpha f
  \|_2\, ,\qquad  \forall f \in C_c^\infty(G)~~.
 \]
\end{lemma}

For the derivation of oscillation estimates, we will need a
local version of this lemma.
\begin{lemma} \label{lem:sobolev}
 Let $K \subset G$ be relatively compact, and suppose that $U \supset \overline{K}$ is open and relatively compact.
There exists a constant $C_{K,U}$ such that for all $f \in
C^\infty(G)$
\begin{equation} \label{eqn:sobolev}
  \| f \|_{\infty, K} \le \binom{2n}{n}^{1/2}  ~ C_{K,U}\big( \sum_{|\alpha| \le n} \|
  X^\alpha f \|_{2,U}^2\big)^{1/2}.
 \end{equation}
\end{lemma}

\begin{prf}
Fix $\psi \in C^\infty  _c(G)$ with $\psi|_K \equiv 1$ and $\psi|_{G
  \setminus U} \equiv 0$.
If $f \in C^\infty(G)$, the Sobolev estimate of the previous lemma
implies the estimate
\begin{eqnarray*}
\| f \|_{\infty, K} & = & \| f \cdot \psi \|_{\infty,K} \\
& \le & C_K \sum_{|\alpha| \le n} \| X^\alpha (\psi \cdot f) \|_{2,G} \\
& = & C_K \sum_{|\alpha| \le n} \| X^\alpha (\psi \cdot f) \|_{2,U} \\
& \le & C_\psi C_K \sum_{|\alpha| \le n} \| X^\alpha f \|_{2,U}~~.
\end{eqnarray*}
Applying the (discrete) Cauchy-Schwarz inequality, we find that
\begin{equation} \label{eqn:sobolev_squareda}
  \| f \|_{\infty, K} \le C_n  ~ C_{K,U}\big( \sum_{|\alpha| \le n} \|
  X^\alpha f \|_{2,U}^2\big)^{1/2}.
 \end{equation}
Here $C_{K,U} = C_K C_\psi$, $C_n =  \mathrm{card}\, \{ \alpha \in
\NN _0^n:
  |\alpha | \leq n \}^{1/2} = \binom{2n}{n}^{1/2}\leq 2^{n-1/2}$, and
  the latter constant depends only on the dimension of $G$.
\end{prf}

In the following it is understood that $C_{K,U}$ is the optimal
constant. The next lemma investigates the behavior of this
constant under translations and dilations.
\begin{lemma} \label{lem:cku}
 Let $K \subset U$, with $K$ and $U$  relatively compact and  $U$
 open. Then

(i)   $C_{xK,xU} = C_{K,U}$ for all  $x\in G$.

(ii) If $0 < r < 1$, then
 \[
 C_{\delta_r(K), \delta_r(U)} \le r^{-Q/2} C_{K,U}~~.
 \]
\end{lemma}
\begin{prf}
(i) is clear,  since  the differential operators $X^\alpha$ are
left-invariant.

(ii) The inequality follows from
\begin{eqnarray*}
 \| f \|_{\infty, \delta_r(K)} & = & \| f \circ \delta_r \|_{\infty,
   K} \\
 & \le & \bbin \,  C_{K,U} \Big(\sum_{|\alpha| \le n} \| X^\alpha (f \circ
 \delta_r)  \|_{2,U}^2\Big)^{1/2} \\
 & = & \bbin C_{K,U} \Big(\sum_{ |\alpha| \le n} r^{2d(\alpha)}\| (X^\alpha f)
 \circ \delta_r \|_{2,U}^2\Big)^{1/2}~~,
\end{eqnarray*}
where we used the homogeneity relation (\ref{eqn:homogen}). As $r
\le 1$, we continue
\begin{eqnarray*} \lefteqn{\bbin C_{K,U} \Big(\sum_{|\alpha| \le n}
  r^{2d(\alpha)}\|   (X^\alpha f) \circ \delta_r
  \|_{2,U}^2\Big)^{1/2}\le} \\
 & \le & \bbin  C_{K,U} \Big(\sum_{ |\alpha| \le n} \| (X^\alpha f)
 \circ \delta_r \|_{2,U}^2 \Big)^{1/2}\\
 & = & \bbin  C_{K,U} r^{-Q/2} \Big(\sum_{ |\alpha| \le n} \|
 (X^\alpha f)  \|_{2,\delta_r(U)} ^2 \Big) ^{1/2}
 ~~,
\end{eqnarray*}
where the last equality is an application of (\ref{eqn:dil_L2}).
Since $ C_{\delta_r(K), \delta_r(U)}$ is the optimal constant, the
conclusion follows.
\end{prf}

The following lemma contains the central estimate for the
oscillation of functions in the Paley-Wiener space. We will write
$B_\epsilon = \{ x \in G : |x|< \epsilon \}$ for the homogeneous
ball of radius $\epsilon $ centered at $e$. All oscillation
estimates below are formulated with respect to these balls, which
are symmetric and
relatively compact sets \cite[Lemma 1.2]{Fo75}. 

\begin{lemma} \label{lem:osc}
There exists a constant $C_G$ depending only on the group $G$ such that for all $f \in E_1$
and all $0<r\leq 1$ the oscillation estimate
 \begin{equation}
\label{eqn:osc} \| {\rm osc}_{B_r} (f) \|_2 \le r C_G \| f \|_2
\end{equation}
holds.
 The constant $C_G>0$ can be chosen to be
\begin{equation} \label{eqn:defCG}
 C_G = \binom{2n}{n}^{1/2}~ 2^{Q/2} b^{Q/2} C_{B_b,B_{2b}} |B_1|^{1/2} \sum_{1 \le |\alpha| \le
 n+1} \| X^\alpha \|_{E_1 \to {\rm L}^2(G)}\, ,
\end{equation}
where $b$ is the constant from Theorem \ref{thm:mean_value}.
\end{lemma}

\begin{prf}
 We first apply Lemma~\ref{thm:mean_value} and estimate
 \begin{eqnarray*}
 \int_G |{\rm osc}_{B_r} f(x) |^2 dx & = & \int_G \sup_{|y|<r} |f(x) -f(xy^{-1})|^2 dx \\
 & \le &  \int_G {\rm sup}_{|y| < r} |y|^2
  \sup_{|z| \le b |y|,1\le j \le n} |X_j f(xz)|^2 dx \\
  & \le & r^2 \int_G \sum_{1\le j \le n} \| X_j f \|_{\infty, xB_{rb}}^2 dx~~,
 \end{eqnarray*}
Next the Sobolev estimate  (\ref{eqn:sobolev}) allows us to
continue
 \begin{eqnarray*}
\ldots & \le & r^2 ~ \binom{2n}{n} \,  \int_G
C_{xB_{rb},xB_{2rb}}^2
\sum_{1\le |\alpha| \le n+1} \| X^\alpha f \|_{2, xB_{2rb}}^2 dx \\
& = & r^2 ~ \binom{2n}{n} ~ C_{B_{rb},B_{2rb}}^2 \int_G \sum_{1\le
|\alpha| \le n+1} \| X^\alpha f \|_{2, xB_{2rb}}^2 dx
~~,\end{eqnarray*} Applying the Lemma~\ref{lem:cku}, we can
continue by
\begin{eqnarray*} \ldots
& \le & r^{2-Q} ~ \binom{2n}{n}  ~ C_{B_{b},B_{2b}}^2 \sum_{1\le
  |\alpha| \le n+1} \int_G \int_{B_{2rb}} |X^{\alpha} f(xy)|^2 dy dx
\\
& = & r^{2-Q} ~ \binom{2n}{n} ~ C_{B_{b},B_{2b}}^2 |B_{2rb}|
\sum_{1\le |\alpha| \le n+1} \| X^{\alpha} f\|^2_2 \\
& = & r^{2-Q} ~\binom{2n}{n} ~ (2rb)^Q C_{B_{b},B_{2b}}^2 |B_{1}| \sum_{1\le |\alpha| \le n+1} \| X^{\alpha} f\|^2_2 \\
& \le & r^2 C_G^2 \| f \|^2_2 ~~.
\end{eqnarray*}
In the last step we have  used the boundedness of differential
operators on $E_1$ by Theorem~\ref{thm:bernstein}.
\end{prf}

Now it is easy  to prove a  sampling theorem for Paley-Wiener spaces.
We first give a version for arbitrary nonuniform sets.

\begin{thm}\label{mainstrat} Given a band-width $\omega >0$, choose
  $s<r< \min(C_G\inv
  \omega   ^{-1/2},\omega^{-1/2})$. Then  every $B_s$-separated and  $B_r$-dense
set $\Gamma \subset G$ is a sampling set for $E_\omega$. In particular
the sampling inequality
\begin{equation}
  \label{eq:c4}
\frac{\omega ^{-Q/2}}{|B_r|^2} (1-r \sqrt{\omega } C_G)^2 \,
\|f\|_2 ^2 \leq \sum _{\gamma \in
  \Gamma } |f(\gamma )|^2 \leq  \frac{\omega ^{-Q/2}}{|B_s|^2} (1+r
  \sqrt{\omega}
C_G)^2 \, \|f\|_2^2
\end{equation}
holds for every $f \in E_\omega $.
\end{thm}
\begin{prf} Assume first that $\omega =1$. Then by Lemma~\ref{lem:osc}
  the
  $r$-oscillation of $f \in E_1$ is at most $rC_G$. Choosing $r <
  \min(C_G\inv,1)$, Theorem \ref{thm:osc_sampl} is applicable  and yields the
  sampling inequality~\eqref{eq:c4}.

The extension to arbitrary $\omega$ is obtained by a dilation argument. For $t>0$
let $U_t$ denote the unitary dilation operator $f \mapsto t^{-Q/2}
f \circ \delta_{t^{-1}}$. 

We claim that $U_t(E_\omega) = E_{t^2 \omega}$. To see this, we
note that by~\eqref{eq:c8} and \eqref{eqn:homogen}  the
sub-Laplacian is $2$-homogeneous and thus  satisfies  $U_t
\mathcal{L} = t^{-2} \mathcal{L} U_t$. Since  the spectral measure
is  unique, we conclude that $U_t \Pi_\mathcal{L}(A) =
\Pi_\mathcal{L}(t^{-2} A) U_t$ for any Borel set $A\subseteq G$
and  $t^{-2} A = \{ t^{-2} r : r \in A \}$. As a consequence,

\begin{equation} \label{eq:c9}
\begin{split}
U_t E_\omega  = & U_t \Pi_\mathcal{L}([0,\omega]) ({\rm L}^2(G)) =
\Pi_\mathcal{L}([0,t^{-2} \omega]) U_t ({\rm L}^2(G)) \\  = &
\Pi_\mathcal{L}([0,t^{-2} \omega]) ({\rm L}^2(G)) = E_{t^{-2}
\omega}~~.
\end{split}
\end{equation}

For the reduction of the general case $\omega >0$ to $\omega =1 $
we choose $t= \sqrt{\omega }$. If $\Gamma \subseteq G$ is
$B_s$-separated and $B_r$-dense, then $\delta _{\sqrt{\omega
}}\Gamma $ is $B_{\sqrt{\omega }s}$-separated and $B_{\sqrt{\omega
}r}$-dense. Since $r \sqrt{\omega } < C_G\inv$, $\delta
_{\sqrt{\omega }} \Gamma $ is a sampling set for $E_1$. Now take
an arbitrary $f\in E_\omega $, then by \eqref{eq:c9} $U_{\delta
_{\sqrt{\omega }}}\in E_1$. By the special case $\omega =1$ we
obtain that
\begin{eqnarray*}
  \lefteqn{\frac{1}{|B_{\sqrt{\omega} r}|^2} (1- r\sqrt{\omega}C_G)^2
  \|U _{\sqrt{\omega }}f\|^2_2 \leq} \\
 & \leq & \sum _{\gamma \in \Gamma} |U_{\sqrt{\omega}} f(\delta
_{\sqrt{\omega}}\gamma )|^2 = \omega ^{-Q/2} \, \sum _{\gamma \in
  \Gamma } |f(\gamma )|^2 \\
&\leq & \frac{1}{|B_{\sqrt{\omega }s}|^2} (1+ r\sqrt{\omega}C_G)^2
  \|U _{\sqrt{\omega }}f\|^2_2 
\end{eqnarray*}
Since $U_t$ is unitary, and $|B_r|^2/|B_{\sqrt{\omega} r}|^2 =
\omega^{-Q} = |B_s|^2/|B_{\sqrt{\omega} s}|^2$, we have proved
\eqref{eq:c4} for all $\omega
>0$.
\end{prf}

\begin{rem}
 At first glance, the dilation property of $E_\omega $ is unexpected,
 therefore  it is instructive to formulate Theorem~\ref{mainstrat}
 explicitly on
  $\RR ^n$.  On $\RR ^n$, $E_\omega = \{ f \in L^2 : \mathrm{supp}
  \hat{f} \subseteq B_{\sqrt{\omega }}\}$. If $f\in E_\omega $, then
  $(U_t f)\, \widehat{} \, (\xi )  = t^{n/2} \hat{f} (t\xi )$ and
  $\mathrm{supp} (U_t f)\, \widehat{} \, \subseteq
  B_{\sqrt{\omega}/t^2}$.

For $G = \RR ^n$, Theorem~\ref{mainstrat} contains a celebrated
theorem of Beurling~\cite{beurling66}. He proved (in our normalization) that
if $r \sqrt{\omega} < \pi /2$
and $\Gamma \subseteq \RR ^n$ is separated and $B_r$-dense, then
 $\Gamma $ is a sampling set for $E_\omega $. The beauty of Beurling's
 Theorem is that the density condition is sharp, whereas our condition
 $r \sqrt{\omega} < C_G\inv$ is  weaker. In our experience,
 oscillation estimates do not lead to sharp density results, their
 strength lies in the general applicability.
\end{rem}

Theorem~\ref{mainstrat} treats the general case of nonuniform
sampling in $E_\omega$. In this case the tightness of the sampling
estimate possesses the upper bound $\frac{|B_{\delta
_r}|^2}{|B_{\delta _s}|^2} \,
\frac{(1+r\sqrt{\omega}C_G)^2}{(1-r\sqrt{\omega}C_G)^2}$. Next we
treat the case of uniform sampling, i.e., sampling on
quasi-lattices, where better estimates can be derived.

First we show that quasi-lattices always exist in a simply
connected solvable Lie group. By contrast, a nilpotent Lie group
allows the existence of a discrete cocompact subgroup only if it
has rational structure constants (see e.g.~\cite{CoGr} for this
well-known theorem of Malcev).

\begin{prop} Let $G$ be a simply connected, connected solvable  Lie
  group. Then there exists a quasi-lattice $\Gamma \subset G$.
\end{prop}
\begin{prf}
We proceed by induction over the dimension of $G$.  The
one-dimensional case is obvious  by taking the lattice $\ZZ \subset
\RR$.

For the induction, choose a  (connected)
normal subgroup  $N \lhd G$ of  codimension one. Then we can  write $G
= N \semdir \RR$, where $\RR$ acts
via a suitable homomorphism $\alpha : \RR \to {\rm Aut}(N)$ and
 the group law on $G \equiv N \semdir  \RR$ is
$(n,t)(m,s) = (n \alpha_t(m),t+s)$. By induction hypothesis, there
exists a quasi-lattice $\Gamma_0 \subset N$ and a (relatively compact)
complement  $C_0 \subset N \subset  G$.

We set
\[ \Gamma  = \{ (\alpha_\ell(\gamma ), \ell) : \gamma \in
\Gamma_0, \ell \in \ZZ \} \]
and
\[
C =  \{ (n,t) :~ n \in C _0, t \in [0,1) \}~~.
\]
We claim that $\Gamma $ is a quasi-lattice in $G$ with complement $C$.

Let $(m,s) \in G= N \semdir \RR$ be arbitrary. We can write  $s= \ell
+t $  for \emph{unique} $\ell \in \ZZ$ and $t \in [0,1)$. Likewise,
since $\Gamma _0$ is a quasi-lattice in $N$, there are unique $\gamma
\in \Gamma _0$ and $n \in C_0$ such that $\gamma n = \alpha _{-\ell }
(m)\in N$. Then $(m,s) = (\alpha _\ell
(\gamma n ), \ell +t ) = (\alpha _\ell (\gamma ), \ell ) (n,t)\in
\Gamma C $. Thus $\Gamma C $ is a covering of $G$, and the uniqueness
of the factorization implies that the covering is disjoint, in other
words, $\Gamma $ is a quasi-lattice.
\end{prf}

The final result of the paper is devoted to regular sampling. We
observe that here the tightness of the sampling estimate
approaches the optimum as the sampling density increases.
\begin{thm}
 Let $\Gamma< G$ be a quasi-lattice, and $\omega>0$. Let $C \subset G$
 be a complement of $\Gamma$ satisfying $C
 \subset B_s$ for a suitable $s > 0$. Then $\delta_r(\Gamma)$ is a sampling set for $E_\omega$,
 as soon as $r$ satisfies $r < \min (s^{-1} \omega^{-1/2}, s^{-1} \omega^{-1/2} C_G)$. The tightness
 of the sampling expansion is $\le (1 +
\frac{2rs \omega^{1/2}}{C_G})^2$.
\end{thm}
\begin{prf} We only consider $\omega=1$, the general case follows
just as in the proof of the previous theorem. As $C$ is a
complement of $\Gamma$, $\delta_r(C)$ is a complement of
$\delta_r(\Gamma)$, due to the fact that $\delta_r$ is an
automorphism. Moreover $C \subset
 B_s$ implies ${\rm osc}_{\delta_r(C)} (f) \le {\rm osc}_{B_{rs}}
 (f)$. Hence (\ref{eqn:sampl_equiv}), with $U=W=\delta_r(C)$, yields the desired statement,
 including the tightness estimate.
\end{prf}

\begin{rem}
 We expect that the arguments employed here for the ${\rm L}^2$-case
 should be adaptable to yield Plancherel-Polya-type results for
 the ${\rm L}^p$-setting, for $1 \le p < \infty$.
\end{rem}


\begin{thebibliography}{99}
\bibitem{AF98}
A.~Aldroubi and H.~G. Feichtinger.
\newblock Exact iterative reconstruction algorithm for multivariate irregularly
  sampled functions in spline-like spaces: the ${L}\sp p$-theory.
\newblock {\em Proc. Amer. Math. Soc.}, 126(9):2677--2686, 1998.

\bibitem{ali-gazeau}
S.~T. Ali, J.-P. Antoine, J.-P. Gazeau, and U.~A. Mueller.
\newblock Coherent states and their generalizations: {A} mathematical overview.
\newblock {\em Rev. Math. Phys.}, 7(7):1013--1104, 1995.


 \bibitem{AnCaDeLe}{P. Aniello, G. Cassinelli, E. De Vito and A. Levrero,
 {\em Wavelet transforms and discrete frames associated to semidirect products,}
 J. Math. Phys. {\bf 39} (1998), 3965-3973.}
\bibitem{BeTa}{D. Bernier and K. Taylor, {\em Wavelets from
  square-integrable representations,} SIAM J. Math. Anal. {\bf 27} (1996),
  594-608.}



\bibitem{ben01}
J.~J. Benedetto and P.~J. S.~G. Ferreira, editors.
\newblock {\em Modern sampling theory}.
\newblock Applied and Numerical Harmonic Analysis. Birkh\"auser Boston Inc.,
  Boston, MA, 2001.
\newblock Mathematics and applications.

\bibitem{ben04}
J.~J. Benedetto and A.~I. Zayed, editors.
\newblock {\em Sampling, wavelets, and tomography}.
\newblock Applied and Numerical Harmonic Analysis. Birkh\"auser Boston Inc.,
  Boston, MA, 2004.

\bibitem{beurling66}
A.~Beurling.
\newblock Local harmonic analysis with some applications to differential
  operators.
\newblock In {\em Some Recent Advances in the Basic Sciences, Vol. 1 (Proc.
  Annual Sci. Conf., Belfer Grad. School Sci., Yeshiva Univ., New York,
  1962--1964)}, pages 109--125. Belfer Graduate School of Science, Yeshiva
  Univ., New York, 1966.

\bibitem{Chr03}{O.~Christensen. {\em An introduction to frames and Riesz
bases.} Birkh\"auser, Bosten, 2003.}

\bibitem{CoGr}{L.~Corwin and F.P.~Greenleaf. {\em Representations
of nilpotent Lie groups and their applications. Part 1: Basic
theory and examples.} Cambridge University Press, Cambridge,
1989.}

\bibitem{Da} {I. Daubechies, {\em The wavelet transform, time-frequency
 localization and signal analysis,}
 IEEE Trans. Inform. Theory {\bf 34} (1988), 961-1005.}

\bibitem{Do}{A.H. Dooley, {\em A nonabelian version of the Shannon sampling
theorem,} Siam. J. Math. Anal. {\bf 20} (1989), 624-633.}

\bibitem{DS52}{ R.J.~Duffin and A.C.~Schaeffer,{\em A class of nonharmonic
Fourier series,} Trans. Am. Math. Soc. {\bf 72} (1952), 341-366.}

\bibitem{FeiGr1} {H.G. Feichtinger and K.-H. Gr\"ochenig, {\em A unified
  approach to atomic decompositions through integrable group representations,}
  52-73 in {\em Function spaces and applications,} Hrsg. M. Cwikel et al.,
  Lecture Notes in Mathematics 1302, Springer, Berlin, 1988.}

\bibitem{FeiGr2}{H.G.~Feichtinger and K.~Gr\"ochenig, {\em
Irregular sampling theorems and series expansions of band-limited
functions,}  J. Math. Anal. Appl. 167, No.2, 530-556 (1992).}


\bibitem{FG92aa}
H.~G. Feichtinger and K.~Gr{\"o}chenig.
\newblock Iterative reconstruction of multivariate band-limited functions from
  irregular sampling values.
\newblock {\em SIAM J. Math. Anal.}, 23(1):244--261, 1992.

\bibitem{FePe}{H.G.~Feichtinger and I.~Pesenson, {\em
Recovery of band-limited functions on manifolds by an iterative
algorithm,}  Heil, Christopher (ed.) et al., {\em Wavelets, frames
and operator theory.} Contemporary Mathematics {\bf 345} (2004),
137-152.}

\bibitem{Fo75}{G.B.~Folland, {\em Subelliptic estimates and
function spaces on nilpotent Lie groups,} Ark. Mat. {\bf 13}
(1975), 161-207.}

\bibitem{FoSt}{G.B.~Folland and E.M.~Stein. {\em Hardy Spaces on
Homogeneous Groups.} Princeton University Press, Princeton, 1982.}


\bibitem{FR05}
M.~Fornasier and H.~Rauhut.
\newblock Continuous frames, function spaces, and the discretization problem.
\newblock {\em J. Fourier Anal. Appl.}, 11(3):245--287, 2005.

\bibitem{Fu} {H. F\"uhr, {\em Wavelet frames and admissibility in
higher dimensions,} J. Math. Phys. {\bf 37} (1996), 6353-6366.}

\bibitem{Fu_LN}{H.~F\"uhr. {\em Abstract Harmonic Analysis of
Continuous Wavelet Transforms.} Springer Lecture Notes in
Mathematics {\bf 1863}, Springer Verlag, Heidelberg, 2005.}

\bibitem{Gr}{K.~Gr\"ochenig, {\em Describing functions: Atomic
decompositions versus frames,} Monatsh. Math. {\bf 112} (1991)
1-42.}

\bibitem{kluvanek}
I.~Kluv{\'a}nek.
\newblock Sampling theorem in abstract harmonic analysis.
\newblock {\em Mat.-Fyz. \v Casopis Sloven. Akad. Vied}, 15:43--48, 1965.


\bibitem{Pe98}{I.~Pesenson, {\em Sampling of Paley-Wiener
functions on stratified groups}, J. Fourier Anal. Appl. {\bf 4}
(1998), 271-281. 
}


\bibitem{Pe01}{I.~Pesenson, {\em Sampling of band-limited vectors,}
J. Fourier Anal. Appl. {\bf 7} (2001), 92-100.}


\bibitem{Pes04}
I.~Pesenson.
\newblock Poincar\'e-type inequalities and reconstruction of {P}aley-{W}iener
  functions on manifolds.
\newblock {\em J. Geom. Anal.}, 14(1):101--121, 2004.

\bibitem{seip04}
K.~Seip.
\newblock {\em Interpolation and sampling in spaces of analytic functions},
  volume~33 of {\em University Lecture Series}.
\newblock American Mathematical Society, Providence, RI, 2004.


\end{thebibliography}

\end{document}